\documentstyle[12pt]{article}
\begin{titlepage}

\title{$\Gamma$-Cohomology and the Selberg Zeta Function}

\author{Ulrich Bunke\thanks{Humboldt-Universit\"at zu Berlin, Institut f\"ur
Reine Mathematik (SFB288), Ziegelstr. 13a, Berlin 10099.
E-mail:ubunke@mathematik.hu-berlin.de
} and
Martin
Olbrich\thanks{Humboldt-Universit\"at zu Berlin, Institut f\"ur Reine
Mathematik (SFB288), Ziegelstr. 13a, Berlin 10099.
E-mail:olbrich@mathematik.hu-berlin.de  }
}

\end{titlepage}

\begin{document}

\newcommand{\oh }{{\bf h}}
\newcommand{\om}{{\bf m}}
\newcommand{\kaaa}{{\bf k}}
\newcommand{\paaa}{{\bf p}}
\newcommand{\taaa}{{\bf t}}
\newcommand{\haaa}{{\bf h}}
\newcommand{\R}{{\bf R}}
\newcommand{\Z}{{\bf Z}}
\newcommand{\C}{{\bf C}}
\newcommand{\HR}{H{\bf R}}
\newcommand{\HC}{H{\bf C}}
\newcommand{\HH}{H{\bf H}}
\newcommand{\PR}{P{\bf R}}
\newcommand{\PC}{P{\bf C}}
\newcommand{\PH}{P{\bf H}}
\newcommand{\laaa}{{\bf l}}
\newcommand{\Saaa}{{\bf S}}
\newcommand{\G}{{\bf G}}
\newcommand{\A}{{\bf A}}
\newcommand{\Naaa}{{\bf N}}
\newcommand{\gaaa}{{\bf g}}
\newcommand{\aaaa}{{\bf a}}
\newcommand{\naaa}{{\bf n}}
\newcommand{\M}{{\bf M}}
\newcommand{\K}{{\bf K}}
\newcommand{\Oaaa}{{\cal O}}
\newcommand{\Haaa}{{\bf H}}
\newcommand{\db}{{\bar{\partial}}}
\newcommand{\Paaa}{{\bf P}}
\newcommand{\g}{{\gaaa}}
 \newcommand{\cZ}{{\cal Z}}
\newcommand{\cE}{{\cal E}}
\newcommand{\cEp}{{{\cal E}^\prime}}
\newcommand{\cU}{{\cal U}}
\newcommand{\Hom}{{\rm Hom}}
\newcommand{\End}{{\rm End}}
\newcommand{\rk}{{\rm rank}}
\newcommand{\im}{{\rm im}}
\newcommand{\spann}{{\rm span}}
\newcommand{\symm}{{\rm symm}}
\newcommand{\cF}{{\cal F}}
\newcommand{\Ree}{{\rm Re}}
\newcommand{\Res}{{\rm Res}}
\newcommand{\op}{{\rm Op}}
\newcommand{\degg}{{\rm deg}}
\newcommand{\Waaa}{{\bf W}}
\newcommand{\Ad}{{\rm Ad}}
\newcommand{\ad}{{\rm ad}}
\newcommand{\codim}{{\rm codim}}
\newcommand{\Gr}{{\rm Gr}}
\newcommand{\coker}{{\rm coker}}
\newcommand{\id}{{\rm id}}
\newcommand{\ord}{{\rm ord}}
\newcommand{\nat}{{\bf N}}
\newcommand{\supp}{{\rm supp}}
\newcommand{\spec}{{\rm spec}}
\newcommand{\Ann}{{\rm Ann}}
\newcommand{\aca}{{\aaaa_\C^\ast}}
\newcommand{\cfu}{{\bf c}}
\newcommand{\ck}{{\cal K}}
\newcommand{\tck}{{\tilde{\ck}}}
\newcommand{\tnk}{{\tilde{\ck}_0}}
\newcommand{\ceep}{{{\cal E}(E)^\prime}}
\newcommand{\ncE}{{{}^\naaa\cE}}
\newcommand{\cB}{{\cal B}}
\newcommand{\hc}{{{\cal HC}(\gaaa,K)}}
\newcommand{\vsl}{{V_{\sigma_\lambda}}}
\newcommand{\czg}{{\cZ(\gaaa)}}
\newcommand{\csl}{{\chi_{\sigma,\lambda}}}
\newcommand{\cR}{{\cal R}}

\maketitle

\newtheorem{prop}{Proposition}[section]
\newtheorem{lem}[prop]{Lemma}
\newtheorem{ddd}[prop]{Definition}
\newtheorem{theorem}[prop]{Theorem}
\newtheorem{kor}[prop]{Corollary}
\newtheorem{ass}[prop]{Assumption}
\newtheorem{con}[prop]{Conjecture}

\tableofcontents
\section{Introduction}
Let $G$ be a semisimple Lie group of real rank one, $K\subset G$ be a maximal
compact subgroup and $\gaaa=\kaaa\oplus \paaa$ be the Cartan decomposition of
the Lie algebra $\gaaa$ of $G$. Let $\aaaa\subset \paaa$ be a one-dimensional
subspace and $M\subset K$ be the centralizer of $\aaaa$ in $K$. Fixing
a positive root system of $(\gaaa,\aaaa)$ we have the Iwasawa decomposition
$\gaaa=\kaaa\oplus\aaaa\oplus\naaa$. Let $(\sigma,V_\sigma)\in\hat{M}$ be
an irreducible representation of $M$ and $\Gamma\subset G$ be a discrete
cocompact torsion-free subgroup.
Then there is a Selberg zeta function $Z_\Gamma(s,\sigma)$, $s\in\aca$,
defined as the analytic continuation of the infinite product
$$Z_\Gamma(s,\sigma)=\prod_{[g]\in
C\Gamma,[g]\not=1,n_\Gamma(g)=1}\prod_{k=0}^\infty
\:det\left(1-e^{(-\rho-s)log(a_g)}S^k(Ad(m_ga_g)_\naaa^{-1})\otimes
\sigma(m_g)\right).$$
Here $C\Gamma$ is the set of conjugacy classes in $\Gamma$, $n_\Gamma(g)$
is maximal number $n\in\nat$ such that $g=h^n$ for some $h\in\Gamma$
and $m_g\in M$, $a_g\in A^+$, $n_g\in N$ are such that $g$ is conjugated in $G$
to $m_ga_gn_g$.
$S^k(Ad(m_ga_g)_\naaa^{-1})$ stands for the $k$'th symmetric power
of $Ad(m_ga_g)^{-1}$ restricted to $\naaa$ and $\rho\in\aaaa^\ast$ is defined
by
$\rho(H):=\frac{1}{2}tr(ad(H)_\naaa)$.
The infinite product converges for $\Ree(s) > \rho$. In this generality
the Selberg zeta function was introduced in \cite{fried86}.

The parameters $(\sigma,\lambda)$ also define a principal series
representation $H^{\sigma,\lambda}$ of $G$.
A Banach globalization is given by
$$H^{\sigma,\lambda}=\{f:G\rightarrow
V_\sigma\:|\:f(gman)=a^{\lambda-\rho}\sigma(m)^{-1}f(g),\:\forall man\in MAN,
f_{|K}\in L^2\},$$
where $G$ acts by the left regular representation.
By $H^{\sigma,\lambda}_{-\omega}$ we denote the space of hyperfunction vectors.

Of main interest are the singularities of the Selberg zeta function, i.e.
the poles and zeros. Their relation to the spectrum of elliptic
differential operators on bundles over $\Gamma\backslash G/K$ and its
compact dual is now well understood (see \cite{bunkeolbrich93},
\cite{bunkeolbrich94}, \cite{bunkeolbrich943} and the forthcoming
\cite{bunkeolbrich95}).
Another description of the singularities in terms of $\naaa$-cohomology was
given in \cite{juhl93}.
S. Patterson \cite{patterson93} conjectured   the relationship of the
singularities of Selberg zeta functions
with the $\Gamma$-cohomology of (subspaces of) principal series
representations.
In the present paper we want to prove the following theorem which settles
this conjecture in the cocompact case.
\begin{theorem}\label{umth1}
The cohomology $H^p(\Gamma,H^{\sigma,\lambda}_{-\omega})$ is finite dimensional
for all $p\ge 0$,
\begin{eqnarray}
\chi(\Gamma,H^{\sigma,\lambda}_{-\omega})=\sum_{p=0}^\infty (-1)^p \dim
H^p(\Gamma,H^{\sigma,\lambda}_{-\omega})&=&0,\label{uass1}\\
\chi(\Gamma,\hat{H}^{\sigma,\lambda}_{-\omega})=\sum_{p=0}^\infty (-1)^p \dim
H^p(\Gamma,\hat{H}^{\sigma,\lambda}_{-\omega})&=&0,\label{uass2}
\end{eqnarray}
and the order of $Z_\Gamma(s,\sigma)$ at $s\in\aca$
can be expressed in terms of the group cohomology of $\Gamma$ with coefficients
in $H^{\sigma,\lambda}_{-\omega}$ as follows :
\begin{eqnarray}\ord_{s=\lambda\not=0} Z_\Gamma(s,\sigma)&=&-\sum_{p=0}^\infty
(-1)^p p \dim H^p(\Gamma,H^{\sigma,\lambda}_{-\omega}),\label{pat}\\
\ord_{s=0} Z_\Gamma(s,\sigma)&=&-\sum_{p=0}^\infty (-1)^p p \dim
H^p(\Gamma,\hat{H}^{\sigma,0}_{-\omega}),
\label{pat0}
\end{eqnarray}
where $\hat{H}^{\sigma,\lambda}_{-\omega}$ is a certain non-trivial extension
of $H^{\sigma,\lambda}_{-\omega}$ with itself.
\end{theorem}
In fact, we propose a new method to study the $\Gamma$- and $\naaa$-cohomology
of the canonical globalizations of arbitrary
Harish-Chandra modules $(\pi,V_{\pi,K})\in \hc$ (that is not restricted to the
rank one case).
Recall (see \cite{wallach92} Ch. 11, \cite{schmid85}, \cite{casselman89}) the
sequence of inclusions
$$V_{\pi,K}\subset    V_{\pi,\omega}\subset    V_{\pi,\infty}\subset
V_{\pi,-\infty}\subset      V_{\pi,-\omega}\subset V_{\pi,for},$$
where $K,\omega,\infty,-\infty,-\omega,for$ stand for $K$-finite, analytic,
smooth, distribution, hyperfunction and formal power series vectors of some
Banach globalization of $V_{\pi,K}$.
For $\omega,\infty,-\infty,-\omega$ these are smooth topological $G$-modules
and the inclusions are continuous.
Our main tool is a resolution (called a standard resolution, see Subsection
\ref{mstr}) of $V_{\pi,-\omega}$ by $\Gamma$-acyclic
and $\naaa$-acyclic (Section \ref{us2})
smooth $G$-modules given by the spaces of smooth sections of homogeneous vector
bundles
over $X=G/K$. The differentials of this resolution are $G$-invariant
differential operators.
The $\Gamma$- or $\naaa$-cohomology of $V_{\pi,-\omega}$ is the cohomology of
the subcomplex
of $\Gamma$- or $\naaa$-invariant vectors of the standard resolution.
A rather simple discussion leads to the finite dimensionality and
Poincar\'e duality (Propositions  \ref{haha}, \ref{mg2}, \ref{ug1}, \ref{ug2}).
Moreover we can show that
$H^\ast(\naaa,V_{\pi,\omega})=H^\ast(\naaa,V_{\pi,K})$ (Proposition
\ref{haha}).
In Proposition \ref{up2}
we construct a long exact sequence relating
$H^\ast(\Gamma,H^{\sigma,\lambda}_{-\omega})$
with the groups $H^\ast(\naaa,V_{\pi,-\omega})$ for all $\pi\in \hat{G}$
with $N_\Gamma(\pi)\not=0$ and with the same infinitesimal character as
$H^{\sigma,\lambda}$.
We finish the proof of Theorem \ref{umth1} by comparing
with the description of the singularities of $Z_\Gamma(s,\sigma)$ given by Juhl
\cite{juhl93}.

Originally, Patterson conjectured Theorem \ref{umth1} for the distribution
globalization of the principal series.
In order to study the cohomology of the distribution globalization
$V_{\pi,-\infty}$
one should take the subcomplex of the standard resolution formed by
the sections with at most exponential growth. Unfortunately we do not know
whether this subcomplex is a resolution of $V_{\pi,-\infty}$  since
we are not able to prove exactness. The essential point to show is the
surjectivity of $(\Omega_G-\mu)$ ($\Omega_G$ is the Casimir operator of $G$ and
$\mu\in\C$) on the space of sections
with at most exponential growth. This is known for trivial bundles
\cite{oshimasaburiwakayama88}.
The problem in the general case is that there is still no
topological  Paley-Wiener theorem for bundles (the non-$K$-finite case).
We believe that concerning the $\naaa$- and $\Gamma$-cohomology (for cocompact
$\Gamma$)
there is no difference between the hyperfunction and distribution
globalizations.
But this difference will certainly appear if one tries to approach the more
general conjecture of Patterson for  finite co-volume or even more general
$\Gamma$'s.

We also have a construction of a standard resolution in the higher rank case.
This, examples and more applications will be the topic of another
paper.\newline
{\bf Acknowledgement :} {\it We are grateful to S.Patterson for explaining to
us his conjecture.
We thank A. Juhl for his constant interest in our work and stimulating
discussions.
The authors are supported by the Sonderforschungsbereich 288
"Differentialgeometrie
und Quantenphysik".}
\section{Acyclicity Lemmas}\label{us2}
\subsection{$\naaa$-acyclicity of $\cE$}
Let $E\rightarrow X$ be the homogeneous vector bundle associated to
the finite dimensional representation $(\gamma,V_\gamma)$ of $K$
and let $\cE$
be its space of smooth sections. $\cE$ is a $\naaa$-module
with the action induced from the left regular action of $G$.
\begin{lem}\label{uli1}
We have
$$H^p(\naaa,\cE)=0,\quad \forall p\ge 1.$$
\end{lem}
{\it Proof:}
Using the Iwasawa decomposition $G=NAK$ we obtain
$$\cE=[C^\infty(G)\otimes V_\gamma]^K=C^\infty(N)\otimes [C^\infty(AK)\otimes
V_\gamma]^K$$
as left $\naaa$-modules, where $\naaa$ acts trivially on $[C^\infty(AK)\otimes
V_\gamma]^K$
and the tensor products are topological ones.
The $\naaa$-cohomology complex
$$\dots\stackrel{d}{\rightarrow}
C^\infty(N)\otimes\Lambda^p\naaa^\ast\stackrel{d}{\rightarrow}
C^\infty(N)\otimes\Lambda^{p+1}\naaa^\ast \stackrel{d}{\rightarrow}\dots$$
can be identified with the de Rham complex of $N$. Since $N\cong
\R^{\dim(\naaa)}$ via the
exponential map we have a contraction $\Phi_t:N\rightarrow N $ of $N$ given by
$\Phi_t(exp(n)):=exp(tn)$, $n\in\naaa$, $t\in [0,1]$. The contraction $\Phi_t$
allows us to define a continuous chain contraction
$$h:C^\infty(N)\otimes\Lambda^p\naaa^\ast\rightarrow
C^\infty(N)\otimes\Lambda^{p-1}\naaa^\ast,\quad p\ge 1$$
satisfying $dh+hd=\id$. The contraction $h$ extends to the tensor product with
$[C^\infty(AK)\otimes V_\gamma]^K$ and the lemma follows.
$\Box$\newline

\subsection{$\naaa$-acyclicity of $\cE(B)$}
Let $B=(\Omega_G-\lambda)^l$ for some $\lambda\in\C$, $l\in\nat$
and $\cE(B)=\{f\in \cE\:|\: Bf=0\}$.
\begin{lem}\label{lez}
We have
$$H^p(\naaa,\cE(B))=0,\quad \forall p\ge 1.$$
\end{lem}
{\it Proof:}
We will use the following fact :
An elliptic operator with real analytic coefficients on an analytic
vector bundle over a non-compact manifold
is surjective on the space of smooth sections of that vector bundle.
By  Lemma \ref{uli1}
$$0 \rightarrow \cE(B)\rightarrow \cE\stackrel{B}{\rightarrow} \cE\rightarrow
0$$
is an $\naaa$-acyclic resolution of $\cE(B)$.
Taking $\naaa$-invariants and using the identification
$$\cE=C^\infty(N)\otimes [C^\infty(AK)\otimes V_\gamma]^K=C^\infty(N)\otimes
C^\infty(A)\otimes V_\gamma$$
we obtain the complex
\begin{equation}\label{ueeqq1}0\rightarrow \ncE(B)\rightarrow
C^\infty(A)\otimes V_\gamma\stackrel{{}^\naaa B}{\rightarrow}
C^\infty(A)\otimes V_\gamma\rightarrow 0.\end{equation}
Here ${}^\naaa B$ is the restriction of $B$
to the subspace of $\naaa$-invariant vectors.
It is a second order translation invariant differential operator on $A$.
The complex (\ref{ueeqq1}) is again exact since ${}^\naaa B$ is still elliptic.
The lemma follows.
$\Box$\newline
\subsection{$\naaa$-acyclicity of $\cE^{for}$ and $\cE^{for}(B)$}
Let $\cE^{for}:=\Hom_{\cU(\kaaa)}(\cU(\g),V_\gamma)$
be the space of formal power series sections of $\cE$.
Then $\cE^{for}$ is a $\g$-module and  hence a $\naaa$-module.
Let $\cE^{for}(B):=\{f\in\cE^{for}\:|\:Bf=0\}$.
\begin{lem}\label{wolf}
We have
$$ H^p(\naaa,\cE^{for})=H^p(\naaa,\cE^{for}(B))=0,\quad \forall p\ge 1.$$
\end{lem}
{\it Proof:}
$\cE^{for}$ is isomorphic to the injective $\naaa$-module
$\Hom(\cU(\naaa)\otimes \cU(\aaaa),V_\gamma)$.
The map
$$ B^*:\cU(\g)\otimes_{\cU(\kaaa)} V_{\tilde\gamma} \rightarrow
\cU(\g)\otimes_{\cU(\kaaa)} V_{\tilde\gamma},$$
given by the multiplication by  $B$,
is injective. This can be seen by going to the graded module
$Gr(\cU(\gaaa)\otimes_{\cU(\kaaa)}V_{\tilde{\gamma}})=S(\paaa)\otimes
V_{\tilde{\gamma}}$.
Hence $B$ is surjective on $\cE^{for}$.
Now one can argue as in the proof of Lemma \ref{lez}. $\Box$\newline
\subsection{$\Gamma$-acyclicity of $\cE$}
Let $\Gamma\subset G$ be a discrete subgroup acting properly
on $X$.
\begin{lem}\label{ugam}
We have
$$H^p(\Gamma,\cE)=0,\quad \forall p\ge 1.$$
\end{lem}
{\it Proof:}
For $p\ge 0$ let $C^p:=\{f:\Gamma^{p+1}\rightarrow \cE\}$ and
$\partial:C^p\rightarrow C^{p+1}$ be defined by
$$(\partial f)(\gamma_0,\dots,\gamma_{p+1}):=\sum_{i=0}^{p+1} (-1)^i
f(\gamma_0,\dots,\check{\gamma}_i,\dots,\gamma_{p+1}), \quad f\in C^p.$$
$\Gamma$ acts on $C^p$ by
$$(\gamma f)(\gamma_0,\dots,\gamma_p)= L_\gamma f(\gamma^{-1}
\gamma_0,\dots,\gamma^{-1} \gamma_p),\quad \gamma\in \Gamma, f\in C^p.$$
Then $H^\ast(\Gamma,\cE)$ is the cohomology of the complex of
$\Gamma$-invariant vectors
$$0\rightarrow {}^\Gamma C^0\stackrel{\partial}{\rightarrow} {}^\Gamma
C^1\stackrel{\partial}{\rightarrow} {}^\Gamma
C^2\stackrel{\partial}{\rightarrow}\dots.$$
Since $\Gamma$ acts properly on $X$ there is an open set $U\subset X$ such that
$\{\gamma U\}_{\gamma\in\Gamma}$ is a locally finite covering of $X$.
Moreover, there is a partition of unity $\{\rho_\gamma\}_{\gamma\in\Gamma}$
such that
$\supp(\rho_\gamma)\in\gamma U$ and $L_{\gamma_1}
\rho_\gamma=\rho_{\gamma_1\gamma}$.
Consider a cocycle $f\in {}^\Gamma C^p$, $\partial f=0$, $p\ge 1$.
Let $F\in C^{p-1}$ be defined by
$$F(\gamma_0,\dots,\gamma_{p-1})(x)=\sum_{\gamma\in\Gamma}
f(\gamma_0,\dots,\gamma_{p-1},\gamma)(x)\rho_\gamma(x).$$
The sum is finite for any $x\in X$.
$F$ is $\Gamma$-invariant, $F\in{}^\Gamma C^{p-1}$. In fact
\begin{eqnarray*}
 L_{\tilde{\gamma}} F(\tilde{\gamma}^{-1} \gamma_0,\dots, \tilde{\gamma}^{-1}
\gamma_{p-1})&=&\sum_{\gamma\in\Gamma}  L_{\tilde{\gamma}}
f(\tilde{\gamma}^{-1} \gamma_0,\dots,\tilde{\gamma}^{-1} \gamma_{p-1},\gamma)
 L_{\tilde{\gamma}} \rho_\gamma \\
&=&\sum_{\gamma\in\Gamma} f(\gamma_0,\dots,\gamma_{p-1},\tilde{\gamma} \gamma)
\rho_{\tilde{\gamma}\gamma }\\
 &=&F(\gamma_0,\dots,\gamma_{p-1}).
\end{eqnarray*}
Moreover
\begin{eqnarray*}
(\partial F)(\gamma_0,\dots,\gamma_p)&=&\sum_{i=0}^p(-1)^i
\sum_{\gamma\in\Gamma}
f(\gamma_0,\dots,\check{\gamma}_i,\dots,\gamma_p,\gamma)\rho_\gamma\\
&=&\sum_{\gamma\in\Gamma} \sum_{i=0}^p(-1)^i
f(\gamma_0,\dots,\check{\gamma}_i,\dots,\gamma_p,\gamma)\rho_\gamma\\
&=&\sum_{\gamma\in\Gamma}(-1)^p   f(\gamma_0,\dots,\gamma_p)\rho_\gamma\\
&=&(-1)^p f(\gamma_0,\dots,\gamma_p).
\end{eqnarray*}
Hence $\partial (-1)^p F=f$.
The lemma follows.
$\Box$\newline

\subsection{$\aaaa\oplus\naaa$-acyclicity of $C^{-\omega}(G)$}
Let $C^{-\omega}(G)$ be the hyperfunctions on $G$ (see \cite{schlichtkrull84}).
We consider $C^{-\omega}(G)$ as a $\aaaa\oplus\naaa$-right module with the
action induced by the right regular representation of $G$.
\begin{lem} \label{unac}
We have
$$H^p(\aaaa\oplus \naaa,C^{-\omega}(G))=0,\quad\forall p\ge 1.$$
\end{lem}
{\it Proof:}
Let $\cB$ be the sheaf of hyperfunctions on $G$. It is a sheaf
of right $\aaaa\oplus \naaa$-modules.
Forming the $\aaaa\oplus \naaa$-cohomology complex locally we obtain the
complex of sheaves
\begin{equation}\label{ueq1}
0\rightarrow \cB\stackrel{d}{\rightarrow} \cB\otimes \Lambda^1
(\aaaa\oplus\naaa)^\ast \stackrel{d}{\rightarrow} \cB\otimes
\Lambda^2(\aaaa\oplus\naaa)^\ast \stackrel{d}{\rightarrow}\dots.\end{equation}
We claim that this complex is exact. By the left $G$-invariance it is enough
to show the exactness at $1\in G$. We employ the Iwasawa decomposition
$G=KAN$. Let $U\subset K$ be a small neighborhood of the identity
which can be identified analytically with an open subset $V\subset
\R^{\dim(K)}$.
The $\aaaa\oplus \naaa$-cohomology complex of $C^{-\omega}(AN)$ can be
identified
with the de Rham complex over $\R^{1+\dim(\naaa)}$.
Thus
$$0\rightarrow \cB_{|UAN}\stackrel{d}{\rightarrow} \cB_{|UAN}\otimes \Lambda^1
(\aaaa\oplus\naaa)^\ast \stackrel{d}{\rightarrow} \cB_{|UAN}\otimes
\Lambda^2(\aaaa\oplus\naaa)^\ast \stackrel{d}{\rightarrow}\dots.$$
is isomorphic to the sheaf version of the partial de Rham complex with
hyperfunction coefficients
on $V \times \R^{1+\dim(\naaa)}$. But this complex of sheaves is exact by
Thm. 3.2 in \cite{komatsu73}.
Thus (\ref{ueq1}) is an exact complex of sheaves. The sheaf of hyperfunctions
on $G$ is flabby. Hence (\ref{ueq1}) is an acyclic resolution of sheaves
with respect to the global section functor.
On the one hand the cohomology groups of the complex of global sections
are the sheaf cohomology groups of $\cB$ and they vanish at all degrees
$p\ge 1$ again because of the flabbyness of $\cB$. On the other hand
they are the $\aaaa\oplus\naaa$-cohomology groups of $C^{-\omega}(G)=\cB(G)$.
The lemma follows.
$\Box$\newline

\subsection{$\Gamma$-acyclicity of $C^{-\omega}(G\times_MV_\sigma)$}
Let $\Gamma\subset G$ be a discrete subgroup such that $\Gamma\backslash G$ is
compact.
Let $(\sigma,V_\sigma)\in \hat{M}$.
We consider $C^{-\omega}(G\times_M V_\sigma)$ as a left $\Gamma$-module with
the action
induced from the left regular action of $G$.
\begin{lem}\label{ugwa}
We have
$$H^p(\Gamma,C^{-\omega}(G\times_MV_\sigma))=0 ,\quad\forall p\ge 0.$$
\end{lem}
{\it Proof:}
Let $\cU$ be a finite open cover of $\Gamma\backslash G$
such that $p:G\rightarrow\Gamma\backslash G$ induces analytic diffeomorphisms
of the connected components of $p^{-1}(U)$ with $U$ for
all $U\in\cU$.
Let $\tilde{\cU}$ be the open cover consisting of the connected components
of the lifts of all $U\in\cU$.
Let $C^p$ be the vector space of \v{C}ech-cochains of the sheaf $\cB_\sigma$
of hyperfunction sections of $G\times_MV_\sigma$
with respect to the cover $\tilde{\cU}$. There is a natural $\Gamma$-action
on $C^p$ given by
$$(\gamma f)(U_0,\dots,U_p)=L_\gamma f(\gamma^{-1}U_0,\dots,\gamma^{-1}U_p),$$
$f\in C^p$, $U_i\in\tilde{\cU}$, $U_0\cap\dots\cap U_p\not=\emptyset$,
$\gamma\in \Gamma$.
Note that as a $\Gamma$-module $C^p=\Hom_\C(\C\Gamma,\C)\otimes V^p$ for a
certain vector space $V^p$
with the trivial $\Gamma$-action.
In fact, let $S^p$ be a set of representatives with respect to
$\Gamma$-translation
of non-trivial intersections of $p+1$ elements of $\tilde{\cU}$.
Then we can choose
$V^p:=\prod_{U\in S^p} \cB_\sigma(U)$.

It follows that $C^p$ is $\Gamma$-acyclic, $H^q(\Gamma,C^p)=0$ for  all $p\ge
0$, $q\ge 1$.
The \v{C}ech  complex
\begin{equation}\label{ueq2}0\rightarrow C^0\stackrel{\partial}{\rightarrow}
C^1\stackrel{\partial}{\rightarrow} C^2
\stackrel{\partial}{\rightarrow}\dots\end{equation}
is $\Gamma$-equivariant. Since $\cB_\sigma$ is flabby (\ref{ueq2}) is  exact at
all degrees $p\ge 1$.
On the one hand the cohomology groups of the
complex of $\Gamma$-invariant vectors in (\ref{ueq2}) are isomorphic to
$H^\ast(\Gamma,C^{-\omega}(G\times_M V_\sigma))$.
On the other hand this complex can be identified with the
\v{C}ech complex of the flabby sheaf of hyperfunction section of
$\Gamma\backslash G\times_MV_\sigma$
with respect to the cover $\cU$. The lemma follows.
$\Box$\newline
\newcommand{\cHC}{{{\cal HC}(\gaaa,K)}}

\section{Resolutions of admissible representations}
\subsection{Constructions of differential operators}

Let $(\pi,V_{\pi,K})\in\cHC$. Then $V_{\pi,K}$ decomposes into a direct sum of
joint generalized eigenspaces of $\cZ(\gaaa)$. Hence we may assume without loss
of generality that there exist $\lambda\in \C$ and $l\in\Naaa$ such that
$B:=(\Omega_G-\lambda)^l\in \Ann(V_{\pi,K})$.

Let $W$ be a finite dimensional $K$-stable subspace of $ V_{\tilde\pi,K}$, the
dual of $V_{\pi,K}$ in the category $\cHC$, which generates $V_{\tilde\pi,K}$
as a  $\cU(\gaaa)$-module. Let $E_0\rightarrow X$ be the homogeneous vector
bundle $G\times_K \tilde W$ and $\cE_0$ be the space of its smooth sections.
Using any globalization $ V_\pi $ of $V_{\pi,K}$ (i.e. a representation of $G$
such that $ V_\pi=V_{\pi,K}$) we can define an embedding
$$
i: V_{\pi,K}\hookrightarrow \cE_0\cong [C^\infty(G)\otimes \tilde W]^K
$$
by
$$\langle i(v)(g),w\rangle:= \langle w,\pi(g^{-1})v\rangle,\qquad v\in
V_{\pi,K},w\in W,g\in G.$$
In fact, the closure of $i(V_{\pi,K})$ in $\cE_0$ is contained in $\cE_0(B)$
and constitutes the maximal globalization $V_{\pi,max}$ of $V_{\pi,K}$ in the
sense of Schmid \cite{schmid85} by its very definition. Schmid's theorem
identifies $V_{\pi,max}$ with the hyperfunction vectors
$V_{\pi,-\omega}:=((V_\pi^\prime)_\omega)^\prime$ of any Banach globalization
$V_\pi$ of $V_{\pi,K}$, hence $V_{\pi,-\omega}$ does not depend on the choice
of the globalization $V_\pi$.

We will also consider the space $V_{\pi,for}:=V_{\tilde\pi, K}^*$ of formal
power series vectors
of $V_{\pi,K}$. There is an exact functor from $\cHC$ to the category of (not
necessarily $K$-finite) $(\gaaa,K)$-modules which sends $V_{\pi,K}$ to
$V_{\pi,for}$. Note that $V_{\pi,for}=\prod_{\gamma\in\hat K}
V_{\pi,K}(\gamma)$.

For homogeneous vector bundles $E$ and $F$ on $X$ we denote by $D(E,F)$ the set
of $G$-invariant differential operators $\cE\rightarrow \cF$.
\begin{prop}\label{1.zeile}
There exist homogeneous vector bundles $E_1,E_2,\dots$ on $X$ and $G$-invariant
differential operators $D_i\in D(E_i,E_{i+1})$, $i=0,1,\dots$, such that the
embedding $i:V_{\pi,-\omega}\hookrightarrow \cE_0(B)$ can be extended to a
(possibly infinite) exact sequence
\begin{equation}\label{*}
 0\rightarrow
 V_{\pi,-\omega}\stackrel{i}{\rightarrow}
\cE_0(B)\stackrel{D_0}{\rightarrow}\cE_1(B)  \stackrel{D_1}{\rightarrow}
\cE_2(B)\stackrel{D_2}{\rightarrow}\dots\ .
\end{equation}
This sequence remains to be exact on the level of formal power series :
\begin{equation}\label{*+}
 0\rightarrow
 V_{\pi,for}\stackrel{i}{\rightarrow}
\cE_0^{for}(B)\stackrel{D_0}{\rightarrow}\cE_1^{for}(B)
\stackrel{D_1}{\rightarrow}
\cE_2^{for}(B)\stackrel{D_2}{\rightarrow}\dots\ .
\end{equation}
\end{prop}
{\it Proof:} Let $\cZ(E)$ be the image of $\cZ(\gaaa)$ in $D(E,E)$. The
following lemma is well known.
\begin{lem}\label{wno}
For any vector bundle $E\rightarrow X$ the $\C[B]$-module $\cZ(E)$ is finitely
generated.
\end{lem}
\begin{lem}\label{isntit}
For any  vector bundle $E\rightarrow X$ we have $\cE(B)_K\in\cHC$.
\end{lem}
{\it Proof:} Let $(\gamma,V_\gamma)$ be the finite dimensional representation
of $K$ corresponding to $E$ and $(\tilde{\gamma},V_{\tilde{\gamma}})$ its dual.
We consider the $K$-equivariant embedding
$$i:V_{\tilde{\gamma}}\hookrightarrow \widetilde{\cE(B)_K}$$
defined by
$$i(\tilde{v})(f):=\langle\tilde{v},f(e)\rangle\ ,$$
where we identify the fibre of $E$ at $e=[K]$
with $V_\gamma$. Let $T:=\cU(\gaaa)(i(V_{\tilde{\gamma}}))$.
For any $t\in T$  the dimension of $\cZ(\gaaa)t$ can be estimated by the
dimension of a generating subspace of the $\C[B]$-module $\cZ(E)$. Thus, by
Lemma \ref{wno}, T is a locally $Z(\gaaa)$-finite and finitely generated
$\cU(\gaaa)$-module. Hence, by a theorem of Harish-Chandra (\cite{wallach88},
3.4.7), $T\in\cHC$. The canonical map $\cE(B)_K\rightarrow \tilde{T}$
is injective by the analyticity of solutions of the equation $Bf=0$. In fact,
an element in the kernel of this map
would have a vanishing Taylor series at $e$.
We obtain that $T\hookrightarrow \widetilde{\cE(B)_K}$ is surjective.
Thus $T=\widetilde{\cE(B)_K}$ and $\cE(B)_K\in\cHC$ since the dual
of a Harish-Chandra module is a Harish-Chandra module, too (\cite{wallach88},
4.3.2).
$\Box$
\begin{lem}\label{hno}
Let $V_{\pi,K}$ be a Harish-Chandra submodule of $\cE(B)_K$. Then there exist a
homogeneous vector bundle $F$ and an operator $D\in D(E,F)$ such that $\ker D
\cap \cE(B)_K=V_{\pi,K}$. We also have $\ker D \cap \cE(B)=V_{\pi,-\omega}$.
\end{lem}
{\it Proof:} According to the proof of Lemma \ref{isntit} there is a surjection
$$ \cU(\gaaa)\otimes_{\cU(\kaaa)} V_{\tilde\gamma}\rightarrow
\widetilde{\cE(B)_K}\ .$$ Let $W$ be a finite dimensional $K$-stable generating
subspace of the Harish-Chandra module
$V_{\pi,K}^\perp\subset\widetilde{\cE(B)_K}$. Then we choose a $K$-equivariant
map $\alpha$ such that the following diagram
$$\begin{array}{ccc}
W&\stackrel{\alpha}{\longrightarrow}&\cU(\gaaa)\otimes_{\cU(\kaaa)}
V_{\tilde\gamma}\\
\downarrow&&\downarrow\\
V_{\pi,K}^\perp&\longrightarrow &\widetilde{\cE(B)_K}
\end{array}$$
commutes. This is possible since
$\cU(\gaaa)\otimes_{\cU(\kaaa)}V_{\tilde\gamma}$ is $K$-semisimple.

We set $F:=G\times_K \tilde W$. The map $\alpha$ can be considered as an
element of
$$[\cU(\gaaa)\otimes_{\cU(\kaaa)} V_{\tilde\gamma}\otimes \tilde W]^K\cong
[\cU(\gaaa)\otimes_{\cU(\kaaa)} \Hom (V_{\gamma},\tilde W)]^K\ .$$
The latter space is canonically isomorphic to $D(E,F)$ via the right regular
representation $R$ of $\cU(\gaaa)$ on $C^\infty(G)\otimes V_\gamma$. Thus
$\alpha$ defines an element $D\in D(E,F)$. If $\alpha(w)=\sum X_i\otimes v_i$,
then
$$ \langle w,Df\rangle_F = \sum \langle v_i, R_{X_i}f\rangle_E \in
C^\infty(G),\ w\in W,v_i\in V_{\tilde\gamma},X_i\in \cU(\gaaa)\ .$$
Let $f\in \cE(B)$, $X\in\cU(\gaaa)$ and $w\in W$. Then we have
\begin{eqnarray}\label{mufti}
\langle w, L_XDf(1)\rangle_F&=&\langle w, DL_Xf\rangle_F\nonumber\\
&=&\sum\langle v_i,R_{X_i}L_Xf(1)\rangle_E\\
&=&\langle w ,L_Xf \rangle_{\cE(B)}\nonumber\\
&=& \langle L_{X^{op}}w ,f \rangle_{\cE(B)}\ ,\nonumber
\end{eqnarray}
where $X\rightarrow X^{op}$ is the anti-automorphism of $\cU(\gaaa)$ induced by
the multiplication with $-1$ on $\gaaa$. By construction  $Df=0$ iff the left
hand side of (\ref{mufti}) vanishes for all $X\in\cU(\gaaa)$ and $w\in W$,
while $f\in V_{\pi,-\omega}$ iff the right hand side does. The lemma
follows.$\Box$\\

In order to prove Proposition \ref{1.zeile} we
iterate Lemma \ref{hno}. $D_i(\cE_i(B)_K)$ is a Harish-Chandra submodule of
$\cE_{i+1}(B)_K$.
Therefore we find a bundle $E_{i+2}$ and an operator $D_{i+1}\in
D(E_{i+1},E_{i+2})$ such that $\ker D_{i+1}\cap\cE_{i+1}(B)_K =
D_i(\cE_i(B)_K)$. We obtain an exact sequence of Harish-Chandra modules
$$ 0\rightarrow
 V_{\pi,K}\stackrel{i}{\rightarrow}
\cE_0(B)_K\stackrel{D_0}{\rightarrow}\cE_1(B)_K
 \stackrel{D_1}{\rightarrow}\cE_2(B)_K
\stackrel{D_2}{\rightarrow}\dots\ .$$
Applying the maximal globalization functor which is exact (see \cite{schmid85})
we end up with (\ref{*}). Analogously, we want to obtain (\ref{*+}) by taking
formal power series vectors. This is possible according to the following lemma.
\begin{lem}\label{form}
For any homogeneous vector bundle $E$ we have
$$ (\cE(B)_K)_{for}=\cE^{for}(B)\ .$$
\end{lem}
{\it Proof:} Let $E$ be associated to the $K$-representation $V_\gamma$. We
have seen in the proof of Lemma \ref{isntit} that there is a surjection
$\cU(\gaaa)\otimes_{\cU(\kaaa)}V_{\tilde\gamma}\rightarrow
\widetilde{\cE(B)_K}$. Hence
\begin{equation}\label{fux} \widetilde{\cE(B)_K}\cong
\cU(\gaaa)\otimes_{\cU(\kaaa)}V_{\tilde\gamma}/(\cE(B)_K)^\perp\  ,
\end{equation}
where $(\cE(B)_K)^\perp$ denotes the annihilator of $\cE(B)_K$ in
$\cU(\gaaa)\otimes_{\cU(\kaaa)}V_{\tilde\gamma}$. Of course, $\cE(B)_K^\perp
\supset B\cU(\gaaa)\otimes_{\cU(\kaaa)}V_{\tilde\gamma}$. We claim that
$\cE(B)_K^\perp = B\cU(\gaaa)\otimes_{\cU(\kaaa)}V_{\tilde\gamma}$. Indeed, if
$\cE(B)_K^\perp$ would be  larger than $
B\cU(\gaaa)\otimes_{\cU(\kaaa)}V_{\tilde\gamma}$, we could find an $f\in\cE_K$
which on the one hand annihilates
$B\cU(\gaaa)\otimes_{\cU(\kaaa)}V_{\tilde\gamma}$, that is $f\in\cE(B)_K$, and
on the other hand is non-zero on $\cE(B)_K^\perp$. This is a contradiction.

Now, the lemma follows by dualizing (\ref{fux}). $\Box$

\subsection{The standard resolution}\label{mstr}

The aim of this subsection is to extend (\ref{*}) to an exact sequence of full
section spaces.
\begin{lem}\label{mlem}
Let $E$, $F$ be homogeneous vector bundles on $X$ and $A\in D(E,F)$
such that $A  \cE(B) =0$.
Then $A = H B$ for some $H\in D(E,F)$.
\end{lem}
{\it Proof:} Let $(\gamma_i,V_{\gamma_i})$, $i=1,2$, be the representations of
$K$ in the fibre of the origin of $E$ and $F$, respectively. Since
$\cE^{for}(B)_K\cong\cE(B)_K$ we have $A\cE^{for}(B)=0$. We consider the
annihilator $\cE^{for}(B)^\perp=\cE(B)_K^\perp$ in
$\cU(\gaaa)\otimes_{\cU(\kaaa)}V_{\tilde\gamma_1}$. We have seen in the proof
of Lemma \ref{form} that this space is equal to $
B\cU(\gaaa)\otimes_{\cU(\kaaa)}V_{\tilde\gamma_1}$.

We consider $A$ as an element of
$[\cU(\gaaa)\otimes_{\cU(\kaaa)}V_{\tilde\gamma_1}\otimes V_{\gamma_2}]^K$.
Therefore $A$ can be written as $\sum A_i\otimes v_i$, where
$A_i\in\cE^{for}(B)$ and $v_i\in V_{\gamma_2}$. Hence there exist elements
$X_i\in \cU(\gaaa)\otimes_{\cU(\kaaa)}V_{\tilde\gamma_1}$ such that $A=\sum
BX_i\otimes v_i$. Set
$$ H:=\sum X_i\otimes v_i\in
[\cU(\gaaa)\otimes_{\cU(\kaaa)}V_{\tilde\gamma_1}\otimes V_{\gamma_2}]^K\cong
D(E,F)\ .$$
Then $A=BH=HB$. $\Box$\\

Let $V_{\pi,K},E_i,D_i$ be as in Proposition \ref{1.zeile}.
\begin{prop}\label{stare}
There exist $H_i\in D(E_i,E_{i+2})$, $i\ge 0$, making the following into an
exact complex:
\begin{equation}\label{**}
0\rightarrow V_{\pi,-\omega}\rightarrow \cE_0\stackrel{{\scriptsize
\left(\begin{array}{c}D_0\\B\end{array}\right)}}
{\longrightarrow}\begin{array}{c}\cE_1\\ \oplus\\ \cE_0\end{array}
                                           \stackrel{
\left({\scriptsize\begin{array}{cc}D_1&H_0\\-B&D_0\end{array}}\right)}
{\longrightarrow}\begin{array}{c}\cE_2\\ \oplus\\\cE_1\end{array}
                                           \stackrel{
\left({\scriptsize\begin{array}{cc}D_2&H_1\\B&D_1\end{array}}\right)}
{\longrightarrow}\dots\ .
\end{equation}
\end{prop}
We shall call (\ref{**}) a standard resolution of $V_{\pi,-\omega}$.\\
{\it Proof:} In order to construct the operators $H_i$ we apply Lemma
\ref{mlem} for $A=D_{i+1}D_i$. Since $B:\cE_i\rightarrow \cE_i$ as an elliptic
operator with analytic coefficients is surjective the exactness of (\ref{**})
is easily reduced to the exactness of (\ref{*}). $\Box$
\newcommand{\ee}{{\rm e}}

\section{$\naaa$-cohomology
% of globalizations of Harish-Chandra modules
}
\subsection{Finite dimensionality}

Let $(\pi,V_{\pi,K}) \in\cHC$.
Recall that $H^*(\naaa,V_{\pi,-\omega})$ carries a natural $MA$-module
structure.
For $\mu\in \aaaa^*_\C$ we define the generalized eigenspace
$$H^p(\naaa,V_{\pi,-\omega})_\mu:= \{\eta\in H^p(\naaa,V_{\pi,-\omega})\:|\:
\exists k \:\mbox{such that}\:\  (H-\mu(H))^k\eta=0 \ \forall H\in \aaaa\}.$$
\begin{prop}\label{haha}
{}
\begin{enumerate}
\item The inclusion $V_{\pi,-\omega}\hookrightarrow V_{\pi,for}$
induces an isomorphism
$$H^p(\naaa,V_{\pi,-\omega})\rightarrow H^p(\naaa,V_{\pi,for})\ .$$
\item $\dim H^p(\naaa,V_{\pi,-\omega})=\dim H^p(\naaa,V_{\pi,for})<\infty\ .$
\item  Assume  $B:=\Omega_G-\lambda \in \Ann (V_{\pi,K})$ for some
$\lambda\in\C$. If $\mu\not=-\rho$, then $\aaaa$ acts semisimply on
$H^p(\naaa,V_{\pi,-\omega})_\mu$.
\end{enumerate}
\end{prop}
{\it Proof:}
According to the Lemmas \ref{lez}, \ref{wolf} and Proposition \ref{1.zeile}
$H^p(\naaa,V_{\pi,*})$ for $*=-\omega$, $for$ is isomorphic to the cohomology
of the subcomplex of $\naaa$-invariants of (\ref{*}) and (\ref{*+}),
respectively. This together with the following lemma implies the proposition.
\begin{lem}\label{hihi}
For any homogeneous vector bundle $E\rightarrow X$ associated to $V_\gamma$ we
have
$$ {}^\naaa\cE(B)={}^\naaa\cE^{for}(B) \ .$$
Furthermore, this space is finite dimensional and consists of elements of the
form
$$ f(n\exp (H)k)=\gamma(k^{-1})\sum P_i(H) e^{\lambda_i(H)},\quad P_i\in
S(\aaaa^*)\otimes V_\gamma,\lambda_i\in \aaaa^*_\C\ .$$
If the assumption of Proposition \ref{haha}, 3. is satisfied, then $\deg P_i\le
1$ and $\deg P_i=0$, whenever $\lambda_i\not=\rho$.
\end{lem}
{\it Proof:}
The $\cU(\aaaa)$-module
$${}^\naaa \cE^{for}(B)\cong (\widetilde {\cE(B)_K}/\naaa(\widetilde
{\cE(B)_K}))^*$$
is finite dimensional (see \cite{wallach88}, Ch.4). Therefore it splits into
generalized weight spaces ${}^\naaa \cE^{for}(B)_\mu, \mu\in\aaaa^*_\C$. $f\in
{}^\naaa \cE^{for}(B)_\mu$, considered as a formal power series on $\aaaa$,
satisfies the differential equations
\begin{equation}\label{muuh}(H+\mu(H))^kf=0\quad  \forall
H\in\aaaa\end{equation}
for a certain $k\in \nat$.
The solutions of (\ref{muuh}) have the form
$$P(H)e^{-\mu(H)},\qquad P\in S(\aaaa^*).$$
They extend to smooth $\naaa$-invariant sections in ${}^\naaa\cE(B)$.

We are left with the proof of the last assertion. We use the following formula
for the Casimir operator applied to $\naaa$-invariant sections of $E$:
$$ \Omega_G f (\exp tH)=(\ee^{t\rho(H)} \frac{d^2}{dt^2} \ee^{-t\rho(H)} -
|\rho|^2 + \gamma(\Omega_M)) f(\exp tH),\ \  H\in\aaaa,|H|=1\ ,
$$
where $\Omega_M$ is the Casimir operator of $M$ and $\gamma$ the
$K$-representation defining $E$. If $f$ is of the form $f(\exp
tH)=P(tH)\ee^{-\mu(tH)}$, $\Omega_G f=\lambda f$  implies
$$\frac{d^2}{dt^2}P(tH)-2(\mu+\rho)(H)\frac{d}{dt}P(tH)+(\langle \mu, \mu
\rangle +2\langle \rho, \mu \rangle  + \gamma(\Omega_M) - \lambda )P(tH)=0.$$
 It follows that $\deg P= 0$  if $\mu\not=-\rho$ and $\deg P\le 1$ in the
remaining case. $\Box$\\
\subsection{Poincar\'e duality}
We consider a double complex of Frechet or dual Frechet spaces
\begin{equation}\label{obenu}\begin{array}{ccccccccccc}\label{otto}
&&&&0&&0&&0&&\\
&&&&\downarrow&&\downarrow&&\downarrow&&\\
&&0&\rightarrow&H^0&\stackrel{d}{\longrightarrow}& H^1&\stackrel{d
}{\longrightarrow}& H^2&\stackrel{d}{\longrightarrow}&\dots\\
&&\downarrow&& \downarrow  &&  \downarrow &&\ \downarrow  &&\\
0&\rightarrow&K^0&\rightarrow&L^{0,0}&
\stackrel{d }{\longrightarrow}& L^{0,1}&\stackrel{d }{\longrightarrow}&
L^{0,2}&\stackrel{d }{\longrightarrow}&\\
&& \downarrow {\scriptstyle \partial}&& \downarrow {\scriptstyle \partial}&&
\downarrow {\scriptstyle \partial}&&
 \downarrow {\scriptstyle \partial}&&\\
0&\rightarrow&K^1&\rightarrow&L^{1,0}&\stackrel{d }{\longrightarrow}&
L^{1,1}&\stackrel{d }{\longrightarrow}& L^{1,2}&\stackrel{d
}{\longrightarrow}&\\
&& \downarrow {\scriptstyle \partial }&&
\downarrow {\scriptstyle \partial}&& \downarrow {\scriptstyle \partial}&&
\downarrow {\scriptstyle \partial}&&\\
0&\rightarrow&K^2&\rightarrow&L^{2,0} &\stackrel{d }{\longrightarrow}&
L^{2,1}&\stackrel{d }{\longrightarrow}& L^{2,2}&
\stackrel{d}{\longrightarrow}&\\
&& \downarrow {\scriptstyle \partial} &&
 \downarrow {\scriptstyle \partial}&& \downarrow {\scriptstyle \partial}&&
\downarrow {\scriptstyle \partial}&&\\
&&\vdots&&&&&&&& .
\end{array}
\end{equation}
such that the horizontal and vertical complexes
\begin{eqnarray*}
&&0\rightarrow K^i \stackrel{d }{\longrightarrow} L^{i,0} \stackrel{d
}{\longrightarrow}  L^{i,1} \stackrel{d }{\longrightarrow}  L^{i,2} \stackrel{d
}{\longrightarrow} \dots\\
&&0\rightarrow H^j \rightarrow L^{0,j} \stackrel{\partial }{\longrightarrow}
L^{1,j} \stackrel{\partial }{\longrightarrow}  L^{2,j} \stackrel{\partial
}{\longrightarrow}\dots
\end{eqnarray*}
are exact  and $K^*,H^*$ have the induced topologies as subspaces.
\begin{lem}\label{mcltrans}
If the differential of the complex
 $$0 \rightarrow H^0 \stackrel{d}{\longrightarrow} H^1 \stackrel{d
}{\longrightarrow}  H^2\stackrel{d}{\longrightarrow}\dots$$
has  closed range, then so does the differential of
$$
 0\rightarrow K^0\stackrel{\partial }{\longrightarrow}  K^1 \stackrel{\partial
}{\longrightarrow}  K^2 \stackrel{\partial }{\longrightarrow} \dots\ .
$$
\end{lem}
 Let $(\pi,V_{\pi,K})\in \hc$.
\begin{prop}\label{mg2}
The $\naaa$-cohomology of $V_{\pi,-\omega}$ satisfies Poincar\'e duality
\begin{equation}\label{ukk1}
H^p(\naaa,V_{\pi,-\omega})^*\cong H^{\dim(\naaa)-p}(\naaa,
V_{\tilde{\pi},\omega})\otimes\Lambda^{\dim{\naaa}}\naaa .
\end{equation}
Moreover,
\begin{equation}\label{ukk2}
H^p(\naaa,V_{\pi,\omega})\cong H^p(\naaa,V_{\pi,K}).
\end{equation}
 \end{prop}
{\it Proof:}
Consider a standard resolution of $V_{\pi,-\omega}$.
By Lemma \ref{lez} the complex (\ref{*}) is a $\naaa$-acyclic resolution.
Taking the $\naaa$-cohomology complex
of the complex (\ref{*}) in the vertical direction we obtain a double complex
of the type (\ref{obenu}).
The first vertical line becomes
 \begin{equation}\label{ukl1}0\rightarrow
V_{\pi,-\omega}\stackrel{d_\naaa}{\rightarrow}V_{\pi,-\omega}\otimes \naaa^*
\stackrel{d_\naaa}{\rightarrow}
V_{\pi,-\omega}\otimes \Lambda^2 \naaa^* \stackrel{d_\naaa}{\rightarrow}\dots
\stackrel{d_\naaa}{\rightarrow}
V_{\pi,-\omega}\otimes \Lambda^ {\dim(\naaa)}\naaa^*\rightarrow 0.
\end{equation}
The dual of this complex is  isomorphic as a complex of $MA$-modules to the
$\naaa$-cohomology complex of
$V_{\tilde{\pi},\omega}$ tensored with $\Lambda^{\dim(\naaa)}\naaa$.
Here we employ the topological duality
$(V_{\pi,-\omega})^\prime=V_{\tilde{\pi},\omega}$.
For (\ref{ukk1}) it is enough to show that the differential of (\ref{ukl1}) has
a closed range.
In view of Lemma \ref{mcltrans} we must show that the differential in
\begin{equation}\label{ulo1}
 0\rightarrow {}^\naaa V_{\pi,-\omega}\rightarrow
{}^\naaa\cE_0(B)\stackrel{D_0}{\rightarrow}{}^\naaa \cE_1(B)
\stackrel{D_1}{\rightarrow}{}^\naaa \cE_2(B)\stackrel{D_2}{\rightarrow}\dots
\end{equation}
has  closed range. But by Lemma \ref{hihi} this complex is finite dimensional.
The isomorphism
(\ref{ukk2}) follows from
the algebraic Poincar\'e duality
$$H^p(\naaa,V_{\pi,K})^*\cong H^{\dim(\naaa)-p}(\naaa,
V_{\tilde{\pi},for})\otimes\Lambda^{\dim(\naaa)}\naaa,$$
(\ref{ukk1})
and \ref{haha}(1).
$\Box$\newline
\section{$\Gamma$-cohomology}
\subsection{Finite dimensionality}
Let $(\pi,V_{\pi,K})\in \hc$ and $\Gamma\subset G$ be a discrete torsion free
cocompact subgroup.
\begin{prop}\label{ug1}
We have
$$\dim H^p(\Gamma,V_{\pi,-\omega})<\infty,\quad \forall p\ge 0.$$
\end{prop}
{\it Proof:}
Let
$$0\rightarrow V_{\pi,-\omega}\rightarrow \cE_0\stackrel{\scriptsize
\left(\begin{array}{c}D_0\\B\end{array}\right)}
{\longrightarrow}\begin{array}{c}\cE_1\\ \oplus\\ \cE_0\end{array}
                                          \stackrel{\scriptsize
\left(\begin{array}{cc}D_1&H_0\\-B&D_0\end{array}\right)}
{\longrightarrow}\begin{array}{c}\cE_2\\ \oplus\\ \cE_1\end{array}
                                           \stackrel{
\scriptsize\left(\begin{array}{cc}D_2&H_1\\B&D_1\end{array}\right)}
{\longrightarrow}\dots
$$
be a standard resolution of $V_{\pi,-\omega}$. By Lemma \ref{ugam}
the cohomology of the subcomplex of $\Gamma$-invariant vectors is isomorphic to
$H^\ast(\Gamma,V_{\pi,-\omega})$.

For any homogeneous vector   bundle $E\rightarrow X$ the space of smooth
$\Gamma$-invariant sections
${}^\Gamma\cE$ can be identified with the space of smooth sections
of $\Gamma\backslash E\rightarrow \Gamma\backslash X$.
Since $B$ is elliptic and normal with respect to the canonical $L^2$-structure
on ${}^\Gamma\cE$ and $\Gamma\backslash X$ is compact we can split
${}^\Gamma\cE={}^\Gamma\cE(B)\oplus {}^\Gamma\cE(B)^\perp$.

We do this splitting for all ${}^\Gamma\cE_i$ entering the standard resolution.
We obtain a complex which is a direct sum of an exact complex built
from the ${}^\Gamma\cE_i(B)^\perp$ and a complex of finite dimensional vector
spaces
\begin{equation}\label{ueq6}0\rightarrow {}^\Gamma\cE_0(B)\stackrel{
\scriptsize\left(\begin{array}{c}D_0\\B\end{array}\right)}
{\longrightarrow}\begin{array}{c}{}^\Gamma \cE_1(B)\\ \oplus\\{}^\Gamma
\cE_0(B)\end{array}
                                           \stackrel{\scriptsize
\left(\begin{array}{cc}D_1&H_0\\-B&D_0\end{array}\right)}
{\longrightarrow}\begin{array}{c}{}^\Gamma\cE_2(B)\\ \oplus\\
{}^\Gamma\cE_1(B)\end{array}
                                           \stackrel{\scriptsize
\left(\begin{array}{cc}D_2&H_1\\B&D_1\end{array}\right)} {\longrightarrow}\dots
\end{equation}
The cohomology of the latter complex is isomorphic to
$H^\ast(\Gamma,V_{\pi,-\omega})$.
The proposition follows. $\Box$\newline
\subsection{Poincar\'e duality}
\begin{prop} \label{ug2}
The $\Gamma$-cohomology of $V_{\pi,\pm\omega}$ satisfies Poincar\'e duality
$$H^p(\Gamma,V_{\pi,-\omega})^\ast\cong H^{n-p}(\Gamma,V_{\tilde\pi,\omega}),$$
where $n=\dim(X)$.
\end{prop}
{\it Proof:}
Note that $V_{\pi,-\omega}$ is a Frechet representation of $\Gamma$
and $V_{\pi,\omega}$ is its topological dual.
Since $\Gamma\backslash X$ is an oriented compact manifold we can find
a finite oriented simplicial complex $P$ being homeomorphic to
$\Gamma\backslash X$.
{}From a baricentric subdivision of $P$ we can construct two oriented
simplicial
complexes $K,\tilde{K}$ being homeomorphic to $P$
such that $\tilde{K}$ is dual to $K$ (see \cite{lefschetz49}, Ch VI).
I.e., for any oriented $p$-simplex $\sigma^p\subset K$ there
is an unique oriented $n-p$-simplex $\tilde{\sigma}^{n-p}\subset \tilde{K}$
such that $\sigma^p\cap\tilde{\sigma}^{n-p}$ is a baricenter of a simplex of
$P$
and the algebraic intersection number is $1$.
Note that $\Gamma\backslash X$ has the homotopy type of the classifying space
$B\Gamma$.
The representation $V_{\pi,-\omega}$ gives rise to a local system over $K$.
We form the associated cochain complex
\begin{equation}\label{ueq5}0\rightarrow C^0\stackrel{\partial}{\rightarrow}
C^1 \stackrel{\partial}{\rightarrow}
C^2\stackrel{\partial}{\rightarrow}\dots,\end{equation}
where
$$C^p=\bigoplus_{\sigma^p\in K} V_{\pi,-\omega}.$$
We have identified the space of constant sections of the local system over
$\sigma^p$
with the fibre over $\sigma^p\cap\tilde{\sigma}^{n-p}$, where
$\tilde{\sigma}^{n-p}$ is dual to $\sigma^p$.
It turns out that the topological dual of the complex (\ref{ueq5})
is exactly the cochain complex
$$0\rightarrow \tilde{C}^0\stackrel{\tilde{\partial}}{\rightarrow} \tilde{C}^1
\stackrel{\tilde{\partial}}{\rightarrow}
\tilde{C}^2\stackrel{\tilde{\partial}}{\rightarrow}\dots$$
associated to $\tilde{K}$ and the local system induced by $V_{\pi,\omega}$.
The pairing $C^p\otimes \tilde{C}^{n-p}\rightarrow \C$ is obtained as follows :
The summand $V_{\pi,-\omega}\subset C^p$ corresponding to $\sigma^p$
is paired nontrivially with the summand $V_{\pi,\omega}\subset \tilde{C}^{n-p}$
corresponding to the dual simplex $\tilde{\sigma}^{n-p}$.
In order to show that Poincar\'e duality holds it is enough to show that
the differential $\partial:C^p\rightarrow C^{p+1}$ has
a closed range for all $p\ge 0$.

Let
$$0\rightarrow L^0\stackrel{d}{\rightarrow}
L^1\stackrel{d}{\rightarrow}L^2\stackrel{d}{\rightarrow}\dots$$
be a standard resolution of $V_{\pi,-\omega}$.
Each $L^p$ is a Frechet representation of $\Gamma$ and gives rise to a local
system over $K$ and to a complex
$$0\rightarrow C^{0,p}\stackrel{\partial}{\rightarrow} C^{1,p}
\stackrel{\partial}{\rightarrow}
C^{2,p}\stackrel{\partial}{\rightarrow}\dots.$$
These complexes fit together to a double complex of Frechet spaces, where
$C^{j,p}\rightarrow C^{j,p+1}$ is induced by $(-1)^p d$.
On the one hand the cohomology with respect to $d$ is concentrated in the zero
degree since the standard resolution is exact.
It yields exactly the complex (\ref{ueq5}).
On the other hand the cohomology with respect to $\partial$
is also concentrated in the zero degree since the $L^p$, $p\ge 0$,
are $\Gamma$-acyclic  by Lemma \ref{ugam}.
It gives  the complex (\ref{ueq6}).
The differential of the latter complex
has closed range. Applying Lemma \ref{mcltrans} we conclude that the
differential
of (\ref{ueq5}) has closed range, too.
$\Box$\newline
In a similar manner one can also prove a Poincar\'e duality using the hermitian
dual.
\section{$\Gamma$-cohomology of the principal series}
\subsection{A long exact sequence}
Let $(\sigma,V_\sigma)\in\hat{M}$ and $\lambda\in\aca$.
We define the representation $\sigma_\lambda$ of $MAN$ on
$V_{\sigma_\lambda}:=V_\sigma$
by $\sigma_\lambda(man)=a^{\rho-\lambda}\sigma(m)$.
Consider $ C^{-\omega}(G)\otimes\vsl$ as a right $\aaaa\oplus\naaa$-module.
Here $\aaaa$ acts on both $C^{-\omega}(G)$ and $\vsl$.
The $\aaaa\oplus\naaa$-cohomology complex of $C^{-\omega}(G)\otimes\vsl$
\begin{equation}\label{uer1}
0\rightarrow C^{-\omega}(G)\otimes\vsl \rightarrow
C^{-\omega}(G)\otimes\vsl\otimes \Lambda^1(\aaaa\oplus\naaa)^\ast\rightarrow
C^{-\omega}(G)\otimes\vsl \otimes
\Lambda^2(\aaaa\oplus\naaa)^\ast\rightarrow\dots
\end{equation}
is exact in all degrees $p\ge 1$ by Lemma \ref{unac}. In fact,
as a right $\aaaa\oplus\naaa$-module $C^{-\omega}(G)\otimes\vsl$
can be identified with $C^{-\omega}(G)\otimes V_\sigma$ with the trivial
$\aaaa\oplus\naaa$-action on $V_\sigma$.
Moreover (\ref{uer1}) admits an $M$-action induced from the right regular
action of $M$ on
$C^{-\omega}(G)$ and $\sigma$. Since $M$ is compact, the subcomplex
\begin{equation}\label{uer2}
0\rightarrow C^{-\omega}(G\times_M\vsl) \rightarrow
C^{-\omega}(G\times_M(\vsl\otimes \Lambda^1(\aaaa\oplus\naaa)^\ast))\rightarrow
C^{-\omega}(G\times_M(\vsl\otimes
\Lambda^2(\aaaa\oplus\naaa)^\ast))\rightarrow\dots
\end{equation}
of $M$-invariants is still acyclic in all degrees $p\ge 1$.
The complex (\ref{uer2}) admits a left $G$-action induced from the left regular
action on $C^{-\omega}(G)$.
By Lemma \ref{ugwa} (\ref{uer2}) is a $\Gamma$-acyclic resolution
of its zero'th cohomology
for any cocompact torsion free discrete subgroup $\Gamma\subset G$.
But the zero'th cohomology of (\ref{uer2}) is the space of $MAN$-invariant
hyperfunctions in $C^{-\omega}(G)\otimes \vsl$ and can be identified as a
$G$-module
with the maximal globalization $H^{\sigma,\lambda}_{-\omega}$ of the principal
series.
\begin{kor}
The cohomology of
\begin{eqnarray}\label{uer3}
\lefteqn{ 0\rightarrow C^{-\omega}(\Gamma\backslash G\times_M\vsl)\rightarrow
C^{-\omega}(\Gamma\backslash G\times_M(\vsl\otimes
\Lambda^1(\aaaa\oplus\naaa)^\ast))\rightarrow}\hspace{5cm}\\
 & &\rightarrow C^{-\omega}(\Gamma\backslash G\times_M(\vsl \otimes
\Lambda^2(\aaaa\oplus\naaa)^\ast))\rightarrow\dots\nonumber
\end{eqnarray}
is isomorphic to $H^\ast(\Gamma,H^{\sigma,\lambda}_{-\omega})$.
\end{kor}
Let $\csl$ be the infinitesimal character of
$H^{\sigma,\lambda}_{-\omega}$.
Let $T$ be the finite set of equivalence classes of irreducible
unitary representations of $G$ with infinitesimal character $\csl$.
The right regular representation of $G$ on $L^2(G)$
decomposes as
\begin{equation}\label{uer4}L^2(\Gamma\backslash
G)=\bigoplus^{\mbox{\scriptsize Hilbert}}_{\pi\in\hat{G}} N_\Gamma(\pi) V_\pi,
\end{equation}
where $N_\Gamma(\pi)\in\nat$   and
$N_\Gamma(\pi)V_\pi:=\bigoplus_{i=1}^{N_\Gamma(\pi)} V_\pi$.
Note that $H^\ast(\naaa,V_{\pi,-\omega})$ carries a natural $MA$-module
structure.
Let $H$ be a unit vector in $\aaaa$.
\begin{prop}\label{up2}
There exists a long exact sequence
\begin{eqnarray}\label{uez35}
\lefteqn{\dots\stackrel{h}{\rightarrow} \bigoplus_{\pi\in T} N_\Gamma(\pi)
[H^p(\naaa,V_{\pi,-\omega})
\otimes \vsl]^M\rightarrow
 H^{p+1}(\Gamma,H^{\sigma,\lambda}_{-\omega})
\rightarrow}\hspace{5cm}\\&&\rightarrow \bigoplus_{\pi\in T}
 N_\Gamma(\pi) [H^{p+1}(\naaa,V_{\pi,-\omega}) \otimes\vsl]^M
\stackrel{h}{\rightarrow}\dots,\nonumber
\end{eqnarray}
where $h:H^p(\naaa,V_{\pi,-\omega}) \otimes\vsl\rightarrow
H^p(\naaa,V_{\pi,-\omega}) \otimes\vsl$
is the action of $H$ induced by the $\aaaa$-module structures of the
$\naaa$-cohomology and of $\vsl$.
\end{prop}
{\it Proof:}
{}From (\ref{uer4}) we obtain a decomposition
$$C^{-\omega}(\Gamma\backslash G\times_M(\vsl \otimes
\Lambda^p(\aaaa\oplus\naaa)^\ast))=\bigoplus_{\pi\in T}
N_\Gamma(\pi)[V_{\pi,-\omega}\otimes \vsl
\otimes\Lambda^p(\aaaa\oplus\naaa)^\ast ]^M \oplus R^p,$$
where $R^p=[C^{-\omega}(\Gamma\backslash G)_\perp\otimes \vsl \otimes
\Lambda^p(\aaaa\oplus\naaa)^\ast]^M$ and $C^{-\omega}(\Gamma\backslash
G)_\perp$
is the space of hyperfunction vectors of the subrepresentation $(\oplus_{\pi\in
T}N_\Gamma(\pi)V_\pi)^\perp\subset L^2(\Gamma\backslash G)$
of the right regular representation.
The complex (\ref{uer3}) decomposes into a direct sum of two subcomplexes
\begin{eqnarray}
\lefteqn{0\rightarrow  \bigoplus_{\pi\in T}
N_\Gamma(\pi)[V_{\pi,-\omega}\otimes \vsl]^M\rightarrow\bigoplus_{\pi\in T}
N_\Gamma(\pi)[V_{\pi,-\omega}\otimes
 \vsl\otimes\Lambda^1(\aaaa\oplus\naaa)^\ast]^M
\rightarrow}\hspace{6cm}\label{uez1}
\\
\hspace{-2cm}& &\rightarrow\bigoplus_{\pi\in T}
N_\Gamma(\pi)[V_{\pi,-\omega}\otimes
\vsl\otimes\Lambda^2(\aaaa\oplus\naaa)^\ast]^M\rightarrow\dots \nonumber
\end{eqnarray}
and
\begin{equation}\label{uer7}\cR^.: 0\rightarrow R^0\stackrel{d}{\rightarrow}
R^1\stackrel{d}{\rightarrow} R^2\stackrel{d}{\rightarrow}\dots.
\end{equation}
\begin{lem}\label{ulk}
The complex (\ref{uer7}) is exact.
\end{lem}
{\it Proof:}
The set of infinitesimal characters $\chi_\pi$ of representations
$\pi\in\hat{G}$
with $N_\Gamma(\pi)\not=0$ has no accumulation points.
Thus we find a finite set $A_i\in \cZ(\gaaa)$ and $\epsilon> 0$
such that
for any $\pi\in\hat{G}$ with $N_\Gamma(\pi)\not=0$
and $\chi_\pi\not=\chi_{\sigma,\lambda}$ there is some $i(\pi)$ with
$|\chi_\pi(A_{i(\pi)})-\chi_{\sigma,\lambda}(A_{i(\pi)})|\ge \epsilon$.
The center $\czg$ acts on the complexes (\ref{uer3}), (\ref{uer7}).
The induced action on
$H^\ast(\Gamma,H^{\sigma,\lambda}_{-\omega})$
has the infinitesimal character $\csl$.
Let $[\alpha]\in H^p(\cR^.)$,
$\alpha\in R^p$, $d \alpha=0$ .
Then for $A\in\cZ(\gaaa)$ we have
$A[\alpha]=[A\alpha]=\csl(A)[\alpha]$.
Thus $(A-\csl(A))\alpha=d \beta(A)$ for some
$\beta(A)\in R^{p-1}$.
Note that $R^{p-1}$ has a further decomposition as a  topological direct sum
$$R^{p-1}=\bigoplus^{\mbox{\scriptsize
topological}}_{\pi\in\hat{G},\chi_\pi\not=\csl}
\bigoplus_{j=1}^{N_\Gamma(\pi)}[V_{\pi,-\omega}\otimes \vsl\otimes
\Lambda^{p-1}(\aaaa\oplus\naaa)^\ast]^M.$$
We decompose
 $$\beta(A)=\bigoplus_{\pi\in\hat{G},\chi_\pi\not=\csl}
\bigoplus_{j=1}^{N_\Gamma(\pi)}\beta^j_\pi(A).$$
Define
$$\gamma:= \bigoplus_{\pi\in\hat{G},\chi_\pi\not=\csl}
\bigoplus_{j=1}^{N_\Gamma(\pi)}
 \frac{\beta^j_\pi(A_{i(\pi)})}{\chi_\pi(A_{i(\pi)})-\chi_{\sigma,\lambda}
(A_{i(\pi)})}.$$
Then $\gamma\in R^{p-1}$ and $\alpha=d\gamma$.
The lemma follows.
$\Box$\newline
{}From Lemma \ref{ulk}
follows that $H^\ast(\Gamma,H^{\sigma,\lambda}_{-\omega})$ is isomorphic to the
cohomology
of (\ref{uez1}).
We have $\Lambda^p(\aaaa\oplus \naaa)^\ast=\Lambda^p\naaa^\ast\oplus
\Lambda^{p-1}\naaa^\ast\otimes \aaaa^\ast$.
We identify $\aaaa^\ast\cong \C$ using $H$.
Then (\ref{uez1}) is isomorphic to a finite direct sum of subcomplexes
\begin{eqnarray}
\lefteqn{0\rightarrow [V_{\pi,-\omega}\otimes \vsl]^M\stackrel{\scriptsize
\left(\begin{array}{c}d \\ H \end{array}\right)}
{\longrightarrow}\begin{array}{c}[V_{\pi,-\omega}\otimes
\vsl\otimes\naaa^\ast]^M\\ \oplus\\
{}[V_{\pi,-\omega}\otimes \vsl]^M\end{array}
\stackrel{\scriptsize \left(\begin{array}{cc}d&0\\-H&d\end{array}\right)}
{\longrightarrow}}\hspace{6cm}\nonumber\\
&&\longrightarrow \begin{array}{c}{}[V_{\pi,-\omega}\otimes\vsl\otimes
\Lambda^2\naaa^\ast]^M\\ \oplus\\{}
[V_{\pi,-\omega}\otimes\vsl\otimes\naaa^\ast]^M \end{array}
                            \stackrel{\scriptsize
\left(\begin{array}{cc}d&0\\H&d\end{array}\right)}
{\longrightarrow}\dots.\label{uez2}
\end{eqnarray}
Here any $\pi\in T$ contributes to (\ref{uez1}) with $N_\Gamma(\pi)$ copies of
(\ref{uez2}).
The differential $d$ is the differential of the $\naaa$-cohomology complex.
A standard argument of  homological algebra now gives the long exact sequence
(\ref{uez35}). This finishes the proof of the proposition. $\Box$\newline
$H^\ast(\naaa,V_{\pi,-\omega})$ is a direct sum of generalized $H$-eigenspaces.
By Proposition \ref{haha} for $\lambda\not=0$ the generalized
$\lambda(H)-\rho(H)$-eigenspace
is in fact $\ker(H-\lambda(H)+\rho(H))$.
\begin{kor}\label{uc1}
For $\lambda\not=0$ ,$p\ge 0$ we have
$$H^p(\Gamma,H^{\sigma,\lambda}_{-\omega})\cong\bigoplus_{\pi\in T}
N_\Gamma(\pi) \left[( H^p(\naaa,V_{\pi,-\omega})\oplus
H^{p-1}(\naaa,V_{\pi,-\omega}))\otimes \vsl\right]^{M,H},$$
where $[.]^{M,H}$ stands for $M$-invariant vectors in the kernel of $H$.
\end{kor}
Again by Lemma  \ref{hihi} $H^2$ acts semisimply on
$H^\ast(\naaa,V_{\pi,-\omega})$
for all $\pi$.
We modify the differential of the complex (\ref{uer1}) replacing the action of
$H$ by $H^2$.
Then it remains still exact in  all degrees $p\ge 1$.
Taking $M$-invariants we obtain a corresponding modification of (\ref{uer2}).
Its zero'th
cohomology is as a $G$-module a non-trivial extension
$\hat{H}^{\sigma,\lambda}_{-\omega}$ of
$H^{\sigma,\lambda}_{-\omega}$
with itself :
$$0\rightarrow  H^{\sigma,\lambda}_{-\omega} \rightarrow
\hat{H}^{\sigma,\lambda}_{-\omega}       \rightarrow
H^{\sigma,\lambda}_{-\omega}\rightarrow 0 .$$
Arguing as above we obtain
\begin{prop}
There exists a long exact sequence
\begin{eqnarray}\label{uez36}
\lefteqn{\dots\stackrel{h^2}{\rightarrow} \bigoplus_{\pi\in T} N_\Gamma(\pi)
[H^p(\naaa,V_{\pi,-\omega})
\otimes \vsl]^M\rightarrow} \hspace{2cm}\\&& \rightarrow
H^{p+1}(\Gamma,\hat{H}^{\sigma,\lambda}_{-\omega})\rightarrow \bigoplus_{\pi\in
T} N_\Gamma(\pi) [H^{p+1}(\naaa,V_{\pi,-\omega}) \otimes\vsl]^M
\stackrel{h^2}{\rightarrow}\dots,\nonumber
\end{eqnarray}
where $h^2:H^p(\naaa,V_{\pi,-\omega}) \otimes\vsl\rightarrow
H^p(\naaa,V_{\pi,-\omega}) \otimes\vsl$
is the action of $H^2$ induced by the $\aaaa$-module structures of the
$\naaa$-cohomology and of $\vsl$.
\end{prop}
\begin{kor}\label{uc2}
For $p\ge 0$ we have
$$H^p(\Gamma,\hat{H}^{\sigma,\lambda}_{-\omega})\cong\bigoplus_{\pi\in T}
N_\Gamma(\pi) \left[( H^p(\naaa,V_{\pi,-\omega})\oplus
H^{p-1}(\naaa,V_{\pi,-\omega}))\otimes \vsl\right]^{M,H^2},$$
where $[.]^{M,H^2}$ stands for $M$-invariant vectors in the kernel of $H^2$.
\end{kor}

\subsection{Proof of the Patterson Conjecture}
The assertions (\ref{uass1}) and (\ref{uass2})
follow immediately from \ref{uc1} and \ref{uc2}.

We now recall the description of the singularities of $Z_\Gamma(s,\gamma)$
given in \cite{juhl93}.
Let $w\in N_K(\aaaa)$ represent the non-trivial element of the Weyl group
$W\cong\Z_2$ of $(\gaaa,\aaaa)$.
Then $w(\lambda)=-\lambda$, $\lambda\in\aca$. Moreover if
$(\sigma,V_\sigma)\in\hat{M}$
we let $\sigma^w$ be the representation of $M$ on $V_\sigma$ given by
$\sigma^w(m)=\sigma(wmw^{-1})$. Thus
$(\sigma_\lambda)^w=(\sigma^w)_{-\lambda}$.
For $(\pi,V_{\pi,K})\in \hc$ let
$$\tilde{\chi}(\pi,\sigma,\lambda):=\sum_{p=0}^{\dim(\naaa)}(-1)^p \dim[
H^p(\naaa,V_{\pi,K})\otimes V_{(\sigma_\lambda)^w}]^{M,H^\infty},$$
where $[.]^{M,H^\infty}$ stands for the $M$-invariants in the
generalized $0$-eigenspace of $H$. By Lemma \ref{hihi} it is enough
to take the kernel of $H^2$.
\begin{theorem}[Juhl, \cite{juhl93} Thm. 7.2.1]\label{jjjuh}
\begin{equation}\label{uhj1}\ord_{s=\lambda}
Z_\Gamma(s,\sigma)=(-1)^{\dim(\naaa)} \sum_{\pi\in\hat{G}} N_\Gamma(\pi)
\tilde{\chi}(\pi,\sigma,\lambda).\end{equation}
\end{theorem}
By Proposition  \ref{mg2}  we have for all $(\pi,V_{\pi,K})\in\hc$
$$H^p(\naaa,
 V_{\pi,K})=H^{\dim(\naaa)-p}(\naaa,V_{\pi,-\omega})^\ast\otimes
\Lambda^{\dim(\naaa)}\naaa^\ast.$$
Let $$\chi(\pi,\sigma,\lambda):=\sum_{p=0}^{\dim(\naaa)}
 (-1)^p\dim[H^p(\naaa,V_{\pi,-\omega})\otimes
V_{(\tilde{\sigma}^w)_\lambda}]^{M,H^2}.$$
Then
$\chi(\pi,\sigma,\lambda)=(-1)^{\dim(\naaa)} \tilde{\chi}(\pi,\sigma,\lambda)$.
It is easy to see from Lemma \ref{ulk} that if
$\chi_\pi\not=\chi_{\tilde{\sigma}^w,\lambda}$,
then $N_\Gamma(\pi)\chi(\pi,\sigma,\lambda)=0$.
Thus we can restrict the summation in (\ref{uhj1}) over the finite set
$\tilde{T}:=\{\pi\in\hat{G}\:|\: \chi_\pi=\chi_{\tilde{\sigma}^w,\lambda}\}$.
Set
$$\chi_1(\Gamma,V_{\pi,-\omega}):=\sum_{p=0}^\infty(-1)^p p
\dim(H^p(\Gamma,V_{\pi,-\omega})).$$
{}From Corollary  \ref{uc1} and \ref{uc2} it follows that
\begin{eqnarray}
\sum_{\pi\in
 \tilde{T}}N_\Gamma(\pi)\chi(\pi,\sigma,\lambda)=
 -\chi_1(\Gamma,H^{\tilde{\sigma}^w,\lambda}_{-\omega}),
&&\lambda\not=0\label{umsta}\\
\sum_{\pi\in
 \tilde{T}}N_\Gamma(\pi)\chi(\pi,\sigma,0)=
-\chi_1(\Gamma,\hat{H}^{\tilde{\sigma}^w,0}_{-\omega}).
     \nonumber
\end{eqnarray}
We apply now identities in the $\Gamma$-cohomology of principal series
representations for different parameters $(\sigma,\lambda)$
in order to reduce (\ref{pat}) and (\ref{pat0}) to (\ref{umsta}) and
(\ref{uhj1}).
It is known that the singularities of $Z_\Gamma(s,\sigma)$ are on
$\aaaa^\ast\cup\imath \aaaa^\ast$.

We first discuss the case $\sigma=\sigma^w$.
If $\lambda\in\aaaa^\ast$, then by applying Poincar\'e  duality twice (one
times
with the complex linear dual and then with the hermitian dual)
we get
 $$\chi_1(\Gamma,H^{\tilde{\sigma}^w,\lambda}_{-\omega})=
(-1)^{\dim(\naaa)+1}\chi_1(\Gamma,H^{\sigma,-\lambda}_{\omega})
=\chi_1(\Gamma,H^{\sigma,\lambda}_{-\omega}).$$
By (\ref{uhj1}) we see that (\ref{pat}) holds for $\lambda\in
\aaaa^\ast\setminus\{0\}$.
In a similar fashion we obtain  (\ref{pat0}).
If $\lambda\in\imath \aaaa^\ast\setminus\{0\}$  we have again by Poincar\'e
duality
 $\chi_1(\Gamma,H^{\tilde{\sigma}^w,\lambda}_{-\omega})
=\chi_1(\Gamma,H^{\sigma,-\lambda}_{-\omega})$.
By (\ref{uhj1}) and the  unitary equivalence
$H^{\sigma,\lambda}=H^{\sigma,-\lambda}$
equation (\ref{pat}) also
holds for imaginary $\lambda$.
Now we consider the case $\sigma\not=\sigma^w$.
Then all singularities of $Z_\Gamma(s,\sigma)$ are on $\imath \aaaa^\ast$.
For $\lambda\in\imath\aaaa^\ast\setminus \{0\}$ we have the unitary equivalence
$H^{\sigma,\lambda}=H^{\sigma^w,-\lambda}$.
Using this and again Poincar\'e duality we obtain
 $\chi_1(\Gamma,H^{\tilde{\sigma}^w,\lambda}_{-\omega})=
\chi_1(\Gamma,H^{\sigma,\lambda}_{-\omega})$.
Now (\ref{pat}) follows for $\lambda\not=0$ from (\ref{uhj1}).
If $\lambda=0$ we apply Poincar\'e duality twice to obtain
 $\chi_1(\Gamma,\hat{H}^{\tilde{\sigma}^w,0}_{-\omega})=
\chi_1(\Gamma,\hat{H}^{\sigma^w,0}_{-\omega})$.
Then (\ref{pat0}) follows since
$\ord_{s=0}Z_\Gamma(s,\sigma)=\ord_{s=0}Z_\Gamma(s,\sigma^w)$
(see \cite{bunkeolbrich943}).
This finishes the proof of Theorem \ref{umth1}.$\Box$\newline

\bibliographystyle{plain}

%\bibliography{literatu}

\begin{thebibliography}{10}

\bibitem{bunkeolbrich93}
U.~Bunke and M.Olbrich.
\newblock The wave kernel for the Laplacian on locally symmetric spaces of rank
  one, theta funktions, trace formulas and the Selberg zeta function.
\newblock Preprint, SFB288, Nr.85, 1993, to appear in Ann. Glob. Anal. Geom.
1994.

\bibitem{bunkeolbrich94}
U.~Bunke and M.~Olbrich.
\newblock Theta and zeta functions for locally symmetric spaces of rank one.
\newblock Preprint 288, No. 118, 1994.

\bibitem{bunkeolbrich943}
U.~Bunke and M.~Olbrich.
\newblock Theta and zeta functions for odd-dimensional locally symmetric spaces
  of rank one.
\newblock Preprint 288, No. 133, 1994.

\bibitem{bunkeolbrich95}
U.~Bunke and M.~Olbrich.
\newblock Theta and Selberg zeta funcions.
\newblock Monograph to appear in the Akademie Verlag, 1995.

\bibitem{casselman89}
W.~Casselman.
\newblock Canonical extensions of Harish-Chandra modules.
\newblock {\em Canadian J. Math.}, 41 (1989), 315--438.

\bibitem{fried86}
D.~Fried.
\newblock The zeta functions of Ruelle and Selberg I.
\newblock {\em Ann. scient. \'ec. norm. sup. $4^{e}$ S\'erie}, 19 (1986),
491--517.

\bibitem{juhl93}
A.~Juhl.
\newblock {\em Zeta-Funktionen, Index-Theorie und hyperbolische Dynamik}.
\newblock  Habilitations\-schrift, Humboldt-Universit\"at zu Berlin,
  1993.

\bibitem{komatsu73}
H.~Komatsu.
\newblock Relative cohomology of sheaves of solutions of differential
  equations.
\newblock In {\em Hyperfunctions and Pseudo-Differential Equations, Proc. of a
  Conf.at Katata 1971, LNM 287}, pages 192--261. Springer-Verlag Berlin
  Heidelberg New York, 1973.

\bibitem{lefschetz49}
S.~Lefschetz.
\newblock {\em Introduction to Topology}.
\newblock Princeton University Press, 1949.

\bibitem{oshimasaburiwakayama88}
T.~Oshima, Y.~Saburi, and M.~Wakayama.
\newblock A note on Ehrenpreis' fundamental principle on a symmetric space.
\newblock In {\em Algebraic Analysis Vol. II}, pages 681--697. Academic Press,
  1988.

\bibitem{patterson93}
S.~Patterson.
\newblock Two conjectures on Kleinian groups.
\newblock Talk at Warwick, March 1993.

\bibitem{schlichtkrull84}
H.~Schlichtkrull.
\newblock {\em Hyperfunctions and Harmonic Analysis on Symmetric Spaces}.
\newblock Birkh\"auser, Boston, 1984.

\bibitem{schmid85}
W.~Schmid.
\newblock Boundary value problems for group invariant differential equations.
\newblock {\em Soc. Math. de France, Asterisque, hors s\'erie}, pages 311--321,
  1985.

\bibitem{wallach88}
N.~R. Wallach.
\newblock {\em Real Reductive Groups}.
\newblock Academic Press, 1988.

\bibitem{wallach92}
N.~R. Wallach.
\newblock {\em Real Reductive Groups II}.
\newblock Academic Press, 1992.

\end{thebibliography}

\end{document}